\documentclass[a4paper,12pt]{article}
\usepackage[cp1251]{inputenc}
\usepackage[T2A]{fontenc}
\usepackage[english]{babel}
\usepackage{amsfonts}
\usepackage[usenames]{color}

\date{}
\newtheorem{Th}{Theorem}[section]

\begin{document}
 \renewcommand{\baselinestretch}{0.96}
 \bf
\begin{center}  How to construct wavelets  on local fields of positive characteristic.
\end{center}
\centerline{ Gleb Sergeevich BERDNIKOV, Iuliia Sergeevna KRUSS,}
\centerline{Sergey Fedorovich LUKOMSKII}

\centerline{Department of Mathematic Analysis, Saratov State University, Saratov, Russia.}

\footnotetext[1]{Correspondence: evrointelligent@gmail.com, KrussUS@gmail.com, LukomskiiSF@info.sgu.ru}
\rm
\large

{\bf Abstract:} We present an algorithm  for construction step wavelets  on local fields of positive characteristic.

 \
 \section*{ Introduction }
   In 2004
  H.Jiang, D.Li, and N.Jin \cite{JLJ} introduced the notion of multiresolution analysis (MRA) on local
  fields $ F^{(s)}$
  of positive characteristic $p$,  proved some properties  and
   constructed "Haar MRA"\ and corresponding
  "Haar wavelets". The wavelet theory  developed in \cite{BJ3, BJ4, BJ1, BJ2, LJ}.
   Construction of non-Haar wavelets is the a basic problem in this theory.
  The problem of constructing orthogonal MRA on the field $ F^{(1)}$ is studied in detail in the works \cite{YuF1,YuF3, YuF2,SL,WP,YuFWP}.
  S.F.Lukomskii, A.M.Vodolazov {\cite{SLAV,AVSL}} considered local field $ F^{(s)}$ as a vector space over the finite field $GF(p^s)$ and constructed  non-Haar wavelets. In \cite{SLAV} the authors construct the mask $m^{(\bf 0)}$ and correspondent refinable function $\varphi$ using some tree with zero as a root. In this case wavelets $\Psi=(\psi^{(\bf l)})_{{\bf l}\in GF(p^s)}$ may be found from the equality
  $$
  \hat\psi^{(\bf l)}=m^{(\bf l)}(\chi)\hat\varphi(\chi{\cal  A}^{-1})
  $$
   where
    ${\cal  A}$ is a dilation operator, $m^{(\bf l)}(\chi)=m^{(\bf 0)}(\chi {\bf r}_0^{-\bf l})$, and
    ${\bf r}_k^{\bf l}$ are Rademacher functions. In the article \cite{GBSF}, the concept of $N$-valid tree was introduced   and an algorithm for constructing  the mask $m^{(\bf 0)}$ and correspondent refinable function $\varphi$ was indicated in the field $F^{(1)}$. In the articles \cite{LBK}, \cite{BKL} the mask $m^{(\bf 0)}$ and correspondent refinable function $\varphi$  were constructed using graph which is obtained from $N$-valid tree by adding new arcs. But in this case we cannot define "masks"\  $m^{(\bf l)}(\chi)$ by the equation  $m^{(\bf l)}(\chi)=m^{(\bf 0)}(\chi {\bf r}_0^{-\bf l})$.

   In this article we give an algorithm for construction of "masks"\  $m^{(\bf l)}(\chi)$ in general case.

 \section{Basic concepts}

Let $p$ be a prime number, $s\in \mathbb N$, $GF(p^s)$ -- finite field. Local field $F^{(s)}$ of positive characteristic $p$ is isomorphic (Kovalski-Pontryagin theorem \cite{GGP}) to the set of formal power series
 $$
 a=\sum_{i=k}^{\infty}{\bf a}_it^i,\ k\in \mathbb{Z},\ {\bf a}_i\in GF(p^s).
 $$

 Addition and multiplication in the field $F^{(s)}$ are defined as sum and product of such series. Therefore
  we will consider local field $F^{(s)}$ of positive
  characteristic $p$ as the field of sequences infinite in both directions
  $$
  a=(\dots ,{\bf 0}_{n-1},{\bf a}_n,{\bf a}_{n+1},\dots),\ {\bf a}_j\in GF(p^s)
  $$
which have only finite number of elements ${\bf a}_j$ with negative $j$
  nonequal to zero, and the operations of addition and multiplication are defined by equalities
 $$
 a\dot+b=(({\bf a}_i\dot+ {\bf b}_i))_{i\in \mathbb Z},
 $$
 \begin{equation} \label{eq1.1}
 ab= (\sum_{i,j:i+j=l}({\bf a}_i {\bf b}_j))_{l\in \mathbb Z},
 \end{equation}
where $"\dot+"$ and $"\cdot"$ are respectively addition and multiplication in $GF(p^s)$.
 The norm of the element $a\in F^{(s)}$ is defined by the equality
 $$
 \|a\|=\|(\dots,{\bf 0}_{n-1},{\bf a}_n,{\bf a}_{n+1},\dots)\|=\left(\frac{1}{p^s}\right)^n, \
 \mbox{\rm если}\  {\bf a}_n\ne {\bf 0}.
 $$
 Therefore
 $$
   F^{(s)}_n=\{a=({\bf a}_j)_{j\in \mathbb Z}: {\bf a}_j\in GF(p^s);\ {\bf a}_j={\bf 0},\ \forall j<n
   \}
   $$
is a ball of radius $p^{-ns}$.

   Neighborhoods $F^{(s)}_n$ are compact subgroups of the group
   $F^{(s)+}$. We will denote them as $F^{(s)+}_n$. They have the following properties:

   1)$\dots\subset F^{(s)+}_1\subset F^{(s)+}_0\subset
   F^{(s)+}_{-1}\subset\dots$

   2)$F^{(s)+}_n/ F^{(s)+}_{n+1}\cong GF(p^s)^+$ и $\sharp (F^{(s)+}_n/ F^{(s)+}_{n+1})=p^s$.

It is noted in \cite{SLAV} that the field $F^{(s)}$ can be described as a linear space over $GF(p^s)$. Using this description one may define the multiplication of element $a\in F^{(s)} $ on element $\mbox{\boldmath$\lambda$} \in GF(p^s)$ coordinatewise, i.e. $\mbox{\boldmath$\lambda$} a =(\dots {\bf 0}_{n-1},\mbox{\boldmath$\lambda$} {\bf a}_n,\mbox{\boldmath$\lambda$} {\bf a}_{n+1},\dots)$, and the modulus
    $\mbox{\boldmath$\lambda$} \in GF(p^s)$ can be defined as
$$ |\mbox{\boldmath$\lambda$}|=\left\{
           \begin{array}{ll}
 1,& \mbox{\boldmath$\lambda$} \ne {\bf 0},\\
 0,& \mbox{\boldmath$\lambda$} = {\bf 0}.\\
  \end{array} \right.
$$

 It is also proved there, that the system $g_k\in F_k^{(s)}\setminus F_{k+1}^{(s)}$ is a basis in $F^{(s)}$, i.e. any element $a\in F^{(s)}$ can be represented as:
    \begin{center}
      $a=\sum\limits_{k\in\mathbb{Z}}\mbox{\boldmath$\lambda$}_kg_k,\
      \mbox{\boldmath$\lambda$}_k\in GF(p^s)$. \\
    \end{center}

    From now on we will consider $g_k=(...,{\bf 0}_{k-1},(1^{(0)},0^{(1)},...,0^{(s-1)})_k,{\bf 0}_{k+1},...)$. In this case $\mbox{\boldmath$\lambda$}_k={\bf a}_k$.
 Let us define the sets
 $$
    H_0^{(s)}=\{h\in G: h={\bf a}_{-1}g_{-1}\dot+{\bf a}_{-2}g_{-2}\dot+\dots \dot+ {\bf a}_{-s}g_{-s}
    \},s\in \mathbb N.
  $$
$$
H_0=\{h\in
G:\;h={\bf a}_{-1}g_{-1}\dot+{\bf a}_{-2}g_{-2}\dot+\dots\dot+{\bf a}_{-s}g_{-s},\;s\in\mathbb
N\}.
$$

 The set $H_0$ is the set of shifts in $F^{(s)}$. It is an analogue of the set of nonnegative integers.

 We will denote the collection of all characters of $F^{(s)+}$ as $X$. The set $X$ generates a commutative group with respect to the multiplication of characters: $(\chi*\phi) (a)=\chi(a)\cdot \phi(a)$. Inverse element is defined as $\chi^{-1}(a)=\overline{\chi(a)}$, and the neutral element is $e(a)\equiv1$.

 Following \cite{SLAV} we define characters $r_n$ of the group $F^{(s)+}$ in the following way.
 Let $x=(\dots,{\bf 0}_{k-1},{\bf x}_k,$ ${\bf x}_{k+1},\dots)$, ${\bf x}_j=(x_j^{(0)},x_j^{(1)},\dots,x_j^{(s-1)})\in GF(p^s)$. The element ${\bf x}_j$ can be written in the form ${\bf x}_j=(x_{js+0},x_{js+1},\dots,x_{js+(s-1)})$. In this case
 $$ x=(...,0,x_{ks+0},x_{ks+1},\dots,x_{ks+s-1},x_{(k+1)s+0},x_{(k+1)s+1},\dots,x_{(k+1)s+s-1},\dots)
 $$
 and the collection of all such sequences $x$ is Vilenkin group. Thus the equality
 $r_n(x)=r_{ks+l}(x)=e^{\frac{2\pi i}{p}(x_{ks+l})}$ defines Rademacher function of $F^{(s)+}$ and every character $\chi\in X$ can be described in the following way:

\begin{equation} \label{eq1.2}
 \chi=\prod \limits_{n\in\mathbb Z}{r}_n^{ a_n},\quad a_n=\overline{0,p-1}.
\end{equation}
The equality (\ref{eq1.2}) can be rewritten as
\begin{equation} \label{eq1.3}
  \chi=
  \prod\limits_{k\in\mathbb Z} r_{ks+0}^{a_k^{(0)}}r_{ks+1}^{a_k^{(1)}}\dots r_{ks+s-1}^{a_k^{(s-1)}}
 \end{equation}
 and let us define
   $$
   r_{ks+0}^{a_k^{(0)}}r_{ks+1}^{a_k^{(1)}}\dots r_{ks+s-1}^{a_k^{(s-1)}}={\bf r}_k^{{\bf a}_k}
   $$
   where ${\bf a}_k=(a_k^{(0)},a_k^{(1)},\dots,a_k^{(s-1)})\in GF(p^s)$. Then (\ref{eq1.3}) takes the form
    \begin{equation} \label{eq1.4}
     \chi =\prod_{k\in \mathbb Z}{\bf r}_k^{{\bf a}_k}.
    \end{equation}

    We will refer to ${\bf r}_k^{(1,0,\dots,0)}={\bf r}_k$ as the Rademacher functions.
    By definition we set

    $$
   ({\bf r}_k^{{\bf a}_k})^{{\bf b}_k}={\bf r}_k^{{\bf a}_k{{\bf b}_k}}, \quad
     \chi^{\bf b}=(\prod {\bf r}_k^{{\bf a}_k})^{\bf b}=\prod {\bf r}_k^{{\bf a}_k\bf b},
     \quad {\bf a}_k, {\bf b}_k, {\bf b}\in
  GF(p^s).
          $$
        It follows that if ${\bf x}=((x_k^{(0)},x_k^{(1)},\dots x_k^{(s-1)}))_{k\in \mathbb Z}$ and ${\bf u}=(u^{(0)},u^{(1)},\dots, u^{(s-1)})\in GF(p^s)$ then
  $$
({\bf r}_k^{{\bf u}},{\bf x})=\prod\limits_{l=0}^{s-1}e^{\frac{2\pi
i}{p}u^{(l)}x_k^{(l)}}.
$$
In \cite{SLAV} the following properties of characters are proved

1) ${\bf r}_k^{{{\bf u}\dot+{\bf v}}}={\bf r}_k^{\bf u}{\bf r}_k^{\bf v}$, ${\bf u}, {\bf v}\in GF(p^s)$.

2) $({\bf r}_k^{\bf v},{\bf u}g_j)=1$, $\forall k\ne j$, ${\bf u}, {\bf v}\in GF(p^s)$.

3)  The set of characters of the field $F^{(s)}$ is a linear space
   $(X,\; *,\; \cdot^{GF(p^s)})$
   over the finite field $GF(p^s)$ with multiplication being an inner operation and the power ${\bf u}\in GF(p^s)$being an outer operation.

4) The set of Rademacher functions $({\bf r}_k)$
   is a basis in the space $(X,\; *,\; \cdot^{GF(p^s)})$.

\noindent
The dilation operator
  ${\cal A}$ in local field $F^{(s)}$  is defined as ${\cal
A}x:=\sum_{n=-\infty}^{+\infty}{\bf a}_ng_{n-1}$, where
 $x=\sum_{n=-\infty}^{+\infty}{\bf a}_ng_n\in F^{(s)}$. In the group of characters it is defined as
 $(\chi {\cal A},x)=(\chi, {\cal
 A}x)$.

 \section{Step Wavelets}
  We will consider a case of scaling function $\varphi$, which generates an orthogonal MRA, being a step function. The set of step functions constant on cosets of a subgroup $F_M^{(s)}$
  with the support ${\rm supp}(\varphi)\subset F^{(s)}_{-N}$ will be denoted as
  $\mathfrak D_M(F^{(s)}_{-N})$, $M,N\in \mathbb N$. Similarly, $\mathfrak D_{-N}({F^{(s)}_{M}}^\bot)$ is a set of step functions, constant on the cosets of a subgroup ${F^{(s)}_{-N}}^{\bot}$
  with the support ${\rm supp}(\varphi)\subset {F^{(s)}_M}^\bot$.

  Let $\varphi\in \mathfrak D_M(F^{(s)}_{-N})$
  generate an orthogonal MRA $\{V_n\}$, satisfies the refinement equation
   $\varphi(x)=\sum_{h\in H_0^{(N+1)}}\beta_h\varphi({\cal
   A}x\dot-h)$\ \cite{SLAV}, which we  rewrite in a frequency from

    \begin{equation}                                      \label{eq2.1}
 \hat\varphi(\chi)=m^{(\bf 0)}(\chi)\hat\varphi(\chi{\cal
  A}^{-1}),
 \end{equation}
 where
 $$
 m^{(\bf 0)}(\chi)=\frac{1}{p^s}\sum_{h\in
 H_0^{(N+1)}}\beta_h\overline{(\chi{\cal A}^{-1},h)}
 $$
 is the mask of equation (\ref{eq2.1}). There exist methods for constructing $ m^{(\bf 0)}(\chi)$ and  $\hat\varphi(\chi)$ (see e.g.\cite{BKL}). We want to construct wavelets $\psi^{\bf (l)}, \ {\bf l}\in GF(p^s), {\bf l}\ne {\bf 0}$ from refinable function $\varphi$.
 We will find these wavelets $\psi^{\bf (l)}$ from the equations
  $$
\hat\psi^{\bf (l)}(\chi)=m^{\bf (l)}(\chi)\hat\varphi(\chi{\cal A}^{-1}),
$$
 and will call the functions $m^{\bf (l)}(\chi)$ masks, too. It is evident that $\hat\psi^{\bf (0)}(\chi)=\hat\varphi(\chi)$.
 \begin{Th}\label{th_1}
 Let $m^{(\bf k)}(\chi)\ ({\bf k}\in GF(p^s))$ be a masks that are constant on  the cosets of a subgroup ${F^{(s)}_{-N}}^{\bot}$ and periodic with any period ${\bf r}_1^{{\bf a}_1}{\bf r}_2^{{\bf a}_2}\dots{\bf r}_\nu^{{\bf a}_\nu}$, ${\bf a}_j\in GF(p^s)$, $\nu\in\mathbb N$. Define wavelets $\psi^{(\bf l)}$ by the equations

$$
\hat\psi^{(\bf l)}(\chi)=m^{(\bf l)}(\chi)\hat\varphi(\chi{\cal A}^{-1}),
$$
where $\varphi\in{\mathfrak D}_{M}(F^{(s)}_{-N})$ is a refinable function.
 The shifts system $(\psi^{(\bf l)}({x\dot-h^{({\bf l})}}))$, ${\bf l}\in GF(p^s), h^{(\bf l)}\in H_0$ will be orthonormal iff
            for any ${\bf a}_{-N}\dots{\bf a}_{-1}\in GF(p^s)$
\begin{equation} \label{eq2.2}
\sum\limits_{{\bf a}_0\in GF(p^s)} m^{(\bf k)}(F^{(s)\bot}_{-N} {\bf r}_{-N}^{{\bf a}_{-N}}\dots {\bf r}_0^{{\bf a}_0})
 \overline{m^{(\bf l)}(F^{(s)\bot}_{-N} {\bf r}_{-N}^{{\bf a}_{-N}}\dots {\bf r}_0^{{\bf a}_0})}=
\delta_{{\bf k},{\bf l}}.
\end{equation}
\end{Th}
  {\bf Proof. The sufficiency}.
 Let $\hat\varphi(\chi)\in\mathfrak
 D_{-N}({F^{(s)}_M}^\bot)$
 Consider scalar product $(\varphi(x\dot-g),\psi^{\bf l}(x\dot-h))$, where $g, h\in H_0$.
 $$
 (\varphi(x\dot-g),\psi^{(\bf l)}(x\dot-h))=\int\limits_{F^{(s)}}\varphi(x\dot-g)\overline{\psi^{(\bf l)}(x\dot-h)}d\,\mu(x)=
 $$
 $$
  =\int\limits_X\hat\varphi_{\cdot\dot-g}(\chi)\overline{\hat\psi^{(\bf l)}_{\cdot\dot-h}(\chi)}
   =\int\limits_X\hat\varphi(\chi)\overline{\hat\varphi(\chi
  A^{-1})}\overline{(\chi,g)}(\chi,h)\overline{m^{(\bf l)}(\chi)}d\,\nu(\chi)=
 $$
 $$
  =\int\limits_{F^{(s)\bot}_{M}}|{\hat\varphi(\chi
  A^{-1})}|^2{(\chi,h\dot-g)}m^{(\bf 0)}(\chi)\overline{m^{(\bf l)}(\chi)}d\,\nu(\chi)=
  $$
  $$
  =\Bigl| h\dot-g=\tilde h={\bf h}_{-1}g_{-1}\dot+{\bf h}_{-2}g_{-2}\dot+\dots\Bigr|=
 $$

 $$
 =\sum_{{\bf a}_{-N}\dots,{\bf a}_0,\dots,{\bf a}_{M-1}}\int\limits_{ F^{(s)\bot}_{-N}
  {\bf r}_{-N}^{{\bf a}_{-N}},\dots, {\bf r}_0^{{\bf a}_0},\dots, {\bf r}_{M-1}^{{\bf a}_{M-1}}}
 |\hat\varphi(F^{(s)\bot}_{-N} {\bf r}_{-N}^{{\bf a}_{-N}}\dots  {\bf r}_{M-1}^{{\bf a}_{M-1}}{\cal A}^{-1})|^2(\chi,\tilde h)\,d\nu(\chi)\cdot
 $$
 $$
 \cdot m^{(\bf 0)}(F^{(s)\bot}_{-N} {\bf r}_{-N}^{{\bf a}_{-N}}\dots {\bf r}_0^{{\bf a}_0})
 \overline{m^{(\bf l)}(F^{(s)\bot}_{-N} {\bf r}_{-N}^{{\bf a}_{-N}}\dots {\bf r}_0^{{\bf a}_0})}d\,\nu(\chi)=
 $$
 $$
  =\sum_{{\bf a}_{-N},\dots,{\bf a}_0}m^{(\bf 0)}(F^{(s)\bot}_{-N} {\bf r}_{-N}^{{\bf a}_{-N}}\dots {\bf r}_0^{{\bf a}_0})
 \overline{m^{(\bf l)}(F^{(s)\bot}_{-N} {\bf r}_{-N}^{{\bf a}_{-N}}\dots {\bf r}_0^{{\bf a}_0})}
 $$
\begin{equation}\label{promres}
 \sum_{{\bf a}_{1},{\bf a}_2,\dots,{\bf a}_{M-1}} |\hat\varphi(F^{(s)\bot}_{-N}
 {\bf r}_{-N}^{{\bf a}_{-N+1}} \dots {\bf r}_0^{{\bf a}_1}\dots
 {\bf r}_{M-2}^{{\bf a}_{M-1}})|^2\int\limits_{F^{(s)\bot}_{-N}
 {\bf r}_{-N}^{{\bf a}_{-N}}\dots {\bf r}_0^{{\bf a}_0}\dots
 {\bf r}_{M-1}^{{\bf a}_{M-1}}}(\chi,\tilde h)d\nu(\chi).
 \end{equation}

By the orthonormality criteria for the system of shifts $(\varphi(x\dot-h))$ of the refinable function $\varphi$ $\forall {\bf a}_{-N},\dots,{\bf a}_{0}\in GF(p^s)$ the following equality holds:
 $$
  \sum_{{\bf a}_{1},{\bf a}_2,\dots,{\bf a}_{M-1}}
  |\hat\varphi(F^{(s)\bot}_{-N} {\bf r}_{-N}^{{\bf a}_{-N+1}} \dots
  {\bf r}_0^{{\bf a}_1}\dots
 {\bf r}_{M-2}^{{\bf a}_{M-1}})|^2=1.
 $$
 Consider integral from (\ref{promres})
 $$
 \int\limits_{F^{(s)\bot}_{-N}
 {\bf r}_{-N}^{{\bf a}_{-N}} \dots
  {\bf r}_0^{{\bf a}_0}\dots
 {\bf r}_{M-1}^{{\bf a}_{M-1}}}(\chi,\tilde h)d\nu(\chi)=\frac{1}{p^{sN}}{\bf 1}_{F^{(s)\bot}_{-N}}(\tilde h) {\bf r}_{-N}^{{\bf a}_{-N}}(\tilde h) \dots{\bf r}_{-1}^{{\bf a}_{-1}}(\tilde h)=
 $$
 $$
 =\frac{1}{p^{sN}}{\bf 1}_{F^{(s)\bot}_{-N}}(\tilde h)\prod\limits_{j=-N}^{-1}e^{\frac{2\pi i}{p}(({\bf h}_j,{\bf a}_j))},
 $$
 where $({\bf h}_j,{\bf a}_j)=h_j^{(0)}a_j^{(0)}+\dots+h_j^{(s-1)}a_j^{(s-1)}$ is a scalar product.

Let us introduce the following notation:
 $$
 m^{(\bf 0)}_{{{\bf a}_{-N}}\dots {{\bf a}_0}}=m^{(\bf 0)}(F^{(s)\bot}_{-N} {\bf r}_{-N}^{{\bf a}_{-N}}\dots {\bf r}_0^{{\bf a}_0}),\ \ m^{(\bf l)}_{{{\bf a}_{-N}}\dots {{\bf a}_0}}=m^{(\bf l)}(F^{(s)\bot}_{-N} {\bf r}_{-N}^{{\bf a}_{-N}}\dots {\bf r}_0^{{\bf a}_0}).
 $$

Then we obtain
 $$
 (\varphi(\cdot\dot-g),\psi^{(\bf l)}(\cdot\dot-h))=
 $$
 $$
 =\frac{1}{p^{sN}}{\bf 1}_{F^{(s)\bot}_{-N}}(\tilde h)
  \sum_{{\bf a}_{-N},\dots,{\bf a}_0}m^{(\bf 0)}_{{{\bf a}_{-N}}\dots {{\bf a}_0}}
 \overline{m^{(\bf l)}_{{{\bf a}_{-N}}\dots {{\bf a}_0}}}
 \prod\limits_{j=-N}^{-1}e^{\frac{2\pi i}{p}(({\bf h}_j,{\bf a}_j))}=
 $$
\begin{equation}\label{sys1}
 =\left\{\begin{array}{lll}
 0&{\it if}& \tilde h\notin F^{(s)\bot}_{-N}; \\
 \frac{1}{p^{sN}}
  \sum\limits_{{\bf a}_{-N},\dots,{\bf a}_0} m^{(\bf 0)}_{{{\bf a}_{-N}}\dots {{\bf a}_0}}
 \overline{m^{(\bf l)}_{{{\bf a}_{-N}}\dots {{\bf a}_0}}}&{\it if}&\tilde h= 0;\\
  \frac{1}{p^{sN}}\sum\limits_{{\bf a}_{-N},\dots,{\bf a}_{-1}}\prod\limits_{j=-N}^{-1}e^{\frac{2\pi i}{p}(({\bf h}_j,{\bf a}_j))}
  \sum\limits_{{\bf a}_0} m^{(\bf 0)}_{{{\bf a}_{-N}}\dots {{\bf a}_0}}
 \overline{m^{(\bf l)}_{{{\bf a}_{-N}}\dots {{\bf a}_0}}}&{\it if}& \tilde h\ne 0,\tilde h\in F^{(s)\bot}_{-N}.
  \end{array}\right.
 \end{equation}
  For $(\psi^{(\bf k)}(x\dot-g),\psi^{(\bf l)}(x\dot-h))$ we can derive similar equality:

$$
 (\psi^{(\bf k)}(\cdot\dot-g),\psi^{(\bf l)}(\cdot\dot-h))=
 $$
 $$
 =\frac{1}{p^{sN}}{\bf 1}_{F^{(s)\bot}_{-N}}(\tilde h)
  \sum_{{\bf a}_{-N},\dots,{\bf a}_0}m^{(\bf k)}_{{{\bf a}_{-N}}\dots {{\bf a}_0}}
 \overline{m^{(\bf l)}_{{{\bf a}_{-N}}\dots {{\bf a}_0}}}
 \prod\limits_{j=-N}^{-1}e^{\frac{2\pi i}{p}(({\bf h}_j,{\bf a}_j))}
 $$
 \begin{equation}\label{sys2}
 =\left\{\begin{array}{lll}
 0&{\it if}& \tilde h\notin F^{(s)\bot}_{-N}; \\
 \frac{1}{p^{sN}}
  \sum\limits_{{\bf a}_{-N},\dots,{\bf a}_0} m^{(\bf k)}_{{{\bf a}_{-N}}\dots {{\bf a}_0}}
 \overline{m^{(\bf l)}_{{{\bf a}_{-N}}\dots {{\bf a}_0}}}&{\it if}&\tilde h= 0;\\
  \frac{1}{p^{sN}}\sum\limits_{{\bf a}_{-N},\dots,{\bf a}_{-1}}\prod\limits_{j=-N}^{-1}e^{\frac{2\pi i}{p}(({\bf h}_j,{\bf a}_j))}
  \sum\limits_{{\bf a}_0} m^{(\bf k)}_{{{\bf a}_{-N}}\dots {{\bf a}_0}}
 \overline{m^{(\bf l)}_{{{\bf a}_{-N}}\dots {{\bf a}_0}}}&{\it if}& \tilde h\ne 0,\tilde h\in F^{(s)\bot}_{-N}.
  \end{array}\right.
 \end{equation}

Thus, if masks $m^{(\bf j)}$ for all ${\bf a}_{-N}\dots{\bf a}_{-1}\in GF(p^s)$ satisfy the condition
$$
\sum\limits_{{\bf a}_0} m^{(\bf k)}_{{{\bf a}_{-N}}\dots {{\bf a}_0}}
 \overline{m^{(\bf l)}_{{{\bf a}_{-N}}\dots {{\bf a}_0}}}=
\delta_{{\bf k},{\bf l}},
$$
then the system of shifts $(\psi^{(\bf l)}({x\dot-h^{({\bf l})}}))$, ${\bf l}\in GF(p^s)$  is an orthonormal system.\\
{\bf The necessity.}
Let us fix ${\bf k,l}\in FG(p^s)$ and consider equalities (\ref{sys1}),(\ref{sys2}) as a system of linear equation with unknowns $x^{{\bf k,l}}_{{\bf a}_{-N}\dots{\bf a}_{-1}}=\sum\limits_{{\bf a}_0} m^{\bf (k)}_{{{\bf a}_{-N}}\dots {{\bf a}_0}}
 \overline{m^{\bf (l)}_{{{\bf a}_{-N}}\dots {{\bf a}_0}}}$ and consider the matrix $A$ of this system.

It is obvious that $A$ is a square matrix $p^{sN}\times p^{sN}$. Let us prove that its determinant is nonequal to zero.

Let us start with $N=1$, $s=1$. In this case

\begin{equation}
A=\frac1p\left(\begin{array}{ccccc}
1 & 1 & 1 & \dots & 1\\
1 & e^{\frac{2\pi i}{p}} & e^{\frac{2\pi i}{p}\cdot 2} &  \dots & e^{\frac{2\pi i}{p}\cdot(p-1)}\\
1 & e^{\frac{2\pi i}{p}\cdot 2} & e^{\frac{2\pi i}{p}\cdot 2\cdot 2} &  \dots & e^{\frac{2\pi i}{p}\cdot 2\cdot(p-1)}\\
\vdots & \vdots &\vdots &\vdots &\vdots \\
1 & e^{\frac{2\pi i}{p}\cdot(p-1)} & e^{\frac{2\pi i}{p}\cdot(p-1)\cdot 2} &  \dots & e^{\frac{2\pi i}{p}\cdot(p-1)\cdot(p-1)}
\end{array}\right)=V,
\end{equation}
where $V$ is Vandermonde matrix, which is known to have nonzero determinant.

For the sake of clarity let us consider a case $N=2$, $s=1$. In this case the matrix $A$ may be represented as block matrix
\begin{equation}
A=\frac1p\left(\begin{array}{ccccc}
V & V & V & \dots & V\\
V & e^{\frac{2\pi i}{p}}V & e^{\frac{2\pi i}{p}\cdot 2}V &  \dots & e^{\frac{2\pi i}{p}\cdot(p-1)}V\\
V & e^{\frac{2\pi i}{p}\cdot 2}V & e^{\frac{2\pi i}{p}\cdot 2\cdot 2}V &  \dots & e^{\frac{2\pi i}{p}\cdot 2\cdot(p-1)}V\\
\vdots & \vdots &\vdots &\vdots &\vdots \\
V & e^{\frac{2\pi i}{p}\cdot(p-1)}V & e^{\frac{2\pi i}{p}\cdot(p-1)\cdot 2}V &  \dots & e^{\frac{2\pi i}{p}\cdot(p-1)\cdot(p-1)}V
\end{array}\right)=V\otimes V,
\end{equation}
 where $\otimes$ symbol corresponds to Kronecker product. By the properties of Kronecker product $\det V\otimes V=(\det V)^p(\det V)^p=(\det V)^{2p}\neq 0$. Thus, again matrix $A$ is nonsingular.

For the case of arbitrary $N$, $s=1$ matrix $A$ can be represented as $A=V\otimes V\otimes\dots\otimes V$ $N$ times and will again have nonzero determinant by the properties of Kronecker product.

Similarly, when $N$ and $s$ are both arbitrary $A=V\otimes V\otimes\dots\otimes V$ $sN$ times.
Thus, the system is nonsingular and has a unique solution, which proves the necessity.
$\square$

Theorem \ref{th_1} can be reformulated in the following way: $m^{(\bf k)}(\chi)$ are the masks of corresponding step compactly supported orthonormal wavelets $\psi^{(\bf l)}(\chi)$ if and only if for each ${\bf a}_{-N}\dots {\bf a}_{-1}\in GF{(p^s)}$ matrix $M({{\bf a}_{-N}\dots {\bf a}_{-1}})$ with elements
$$M_{\bf l, a_0}({{\bf a}_{-N},\dots, {\bf a}_{-1}})=m^{\bf (l)}(F^{(s)\bot}_{-N} {\bf r}_{-N}^{{\bf a}_{-N}}\dots {\bf r}_0^{{\bf a}_0})$$
is unitary. The sufficiency of this theorem was proved in \cite{JLJ} (theorem 3). For step refinable functions
the condition (\ref{eq2.2}) is necessary and sufficient. If the condition (\ref{eq2.2}) is fulfilled then the functions  $
\hat\psi^{\bf (l)}(\chi)=m^{\bf (l)}(\chi)\hat\varphi(\chi{\cal A}^{-1})
$ form a wavelet system \cite{JLJ}.
For  a step refinable function we can describe an algorithm for constructing masks $m^{\bf (l)}$ and wavelets $\psi^{\bf (l)}$ , ${\bf l}\in GF{(p^s)}$.

Let us assume we have all the values of $m^{\bf (0)}(\chi)$. We may obtain them using an algorithm presented in \cite{BKL}. Recall the notation:
$$
 m^{\bf (0)}_{{{\bf a}_{-N}}\dots {{\bf a}_0}}=m^{\bf (0)}(F^{(s)\bot}_{-N} {\bf r}_{-N}^{{\bf a}_{-N}}\dots {\bf r}_0^{{\bf a}_0}),\ \ m^{\bf (l)}_{{{\bf a}_{-N}}\dots {{\bf a}_0}}=m^{\bf (l)}(F^{(s)\bot}_{-N} {\bf r}_{-N}^{{\bf a}_{-N}}\dots {\bf r}_0^{{\bf a}_0}).
$$

{\bf 1)} For each ${\bf a}_{-N}\dots {\bf a}_{-1}$ we construct a matrix $M({{\bf a}_{-N}\dots {\bf a}_{-1}})\in Mat _{p^s\times p^s}(\mathbb C)$ with elements $M_{{\bf l},{\bf a}_0}({{\bf a}_{-N}\dots {\bf a}_{-1}})$ the following way.
The first row consists of all the values

 $$
 m^{\bf (0)}_{{{\bf a}_{-N}}\dots {{\bf a}_{-1}}, 0}, m^{\bf (0)}_{{{\bf a}_{-N}}\dots {{\bf a}_{-1}}, 1},
 \dots ,m^{\bf (0)}_{{{\bf a}_{-N}}\dots {{\bf a}_{-1}}, p^s-1}
 $$
    where ${\bf a}_{-N}\dots {\bf a}_{-1}$ are fixed and $j=a_0^{(0)}+a_0^{(1)}p+\dots +a_0^{(s-1)} p^{s-1}$ calculated from ${\bf a_0}=(a_0^{(0)},a_0^{(1)},\dots a_0^{(s-1)}) $. Supplement this matrix to unitary in the following way.

    If $m^{\bf (0)}_{{{\bf a}_{-N}}\dots {{\bf a}_{-1}}, 0}\ne 0$  then we make
    $M_{{\bf l},{\bf l}}=1$ for ${\bf l}\ne {\bf 0}$ and $M_{{\bf l},{\bf a_0}}=0$ for
    ${\bf l}\ne {\bf 0},{\bf l}\ne {\bf a}_0$.

     If $m^{\bf (0)}_{{{\bf a}_{-N}}\dots {{\bf a}_{-1}}, 0}= 0$ then there exists number
    $$
    j=j({\bf a}_{0})=a_0^{(0)}+a_0^{(1)}p+\dots +a_0^{(s-1)} p^{s-1}
    $$
        for which $m^{\bf (0)}_{{{\bf a}_{-N}}\dots {{\bf a}_{-1}}, j}\ne 0$. This nonzero value exists by the property of $m^{\bf (0)}$ (see e.g.\cite{JLJ} ) In this case we  make $ M_{ {\bf j},{\bf 0}}=1$, $M_{{\bf l},{\bf l}}=1$ for ${\bf l}\ne {\bf 0}, {\bf l}\ne {\bf j}$, and $M_{{\bf l},{\bf a_0}}=0$ in another case.

{\bf 2)} Run the Gram-Schmidt process on each matrix in order to make them unitary.

{\bf 3)} Now for each ${\bf l}\in GF(p^s),\ {\bf l}\ne {\bf 0} $ we find the values of the mask $m^{\bf (l)}$ from the equalities
$$m^{\bf (l)}(F^{(s)\bot}_{-N} {\bf r}_{-N}^{{\bf a}_{-N}}\dots {\bf r}_0^{{\bf a}_0})=M_{\bf l, a_0}({{\bf a}_{-N}\dots {\bf a}_{-1}}).
$$.
{\bf 4)} The wavelets $\psi^{\bf (l)}$ can be obtained using the formula
$$
\hat\psi^{\bf (l)}(\chi)=m^{\bf (l)}(\chi)\hat\varphi(\chi{\cal A}^{-1})
$$
and performing inverse Fourier transform.

First and second authors  have performed the work of the state task of
Russian Ministry of Education and Science (project 1.1520.2014K).
The third author was supported RFBR, grant 16-01-00152.

\end{document}